\documentclass[12pt,a4paper]{amsart}%

\usepackage{amsmath}
\usepackage{fullpage}
\usepackage[final]{pdfpages}

\newtheorem{theorem}{Theorem}
\newtheorem{notation}{Notation}

\newcommand{\NextVer}[1]{}

\newcommand{\Sc}{{\mathcal{S}}}

\newtheorem{cor}[theorem]{Corollary}
\newtheorem{proposition}[theorem]{Proposition}

\newtheorem*{defn}{Definition}
\newtheorem*{remark*}{Remark}

\newcommand{\eps}{\varepsilon}

\newcommand{\C}{\mathbb{C}}

\newcommand{\irr}{\mathrm{irr}}

\begin{document}

\author[Aizenbud]{Avraham Aizenbud}
\address{Avraham Aizenbud,
Faculty of Mathematical Sciences,
Weizmann Institute of Science,
Rehovot, Israel}
\email{aizenr@gmail.com}
\urladdr{http://aizenbud.org}

\date{\today}

%
%
%
%
%
%
%
%

\title[Necessary and sufficient Gelfand-Kazhdan criterion]{The Gelfand-Kazhdan criterion as a necessary and sufficient criterion}
\maketitle
\tableofcontents

\begin{abstract}
We show that under certain conditions the Gelfand-Kazhdan criterion for the Gelfand property is a necessary condition.  We work in the generality of finite groups, however part of the argument carries over to $p$-adic and real groups.
\end{abstract}
\section{Introduction}
In this note we study the Gelfand-Kazhdan criterion for the Gelfand property (see \cite{GK}) and show that under some conditions it is not only a sufficient condition but also a necessary one.
We discus mostly finite groups, but we hope that some of these methods can be pushed to the generality of $p$-adic groups (and even Lie groups).
The Gelfand-Kazhdan criterion was originally developed as a version of the Gelfand trick 
that is valid for $p$-adic groups and not only for compact groups. However,  even for finite groups  the Gelfand-Kazhdan criterion is slightly  more informative than the Gelfand trick.

The main result of this note is the following:
\begin{theorem}
\label{thm:GK1}
Let $G$ be a finite group,  $\theta:G\to G$ be an involution and $H \subset G$ be a $\theta$-stable subgroup. Assume that for any $x\in G$ there exists $g\in G$ s.t.
$gx^{-1}g^{-1}=\theta(x)$.
Then the following are equivalent:
\begin{enumerate}
\item $(G,H)$ is a Gelfand pair.
\item  For any $x\in G$ there are $h_1,h_2\in H$ s.t. $h_1x^{-1} h_2 =\theta(x)$.
\end{enumerate}
\end{theorem}

We also have slightly more general version of this theorem.
\begin{theorem}
\label{thm:GK2}
Let $G$ be a finite group,  $\theta:G\to G$ be an involution and $H \subset G$ be a $\theta$-stable subgroup.
Then the following are equivalent:
\begin{enumerate}
\item $(G,H)$ is a Gelfand pair and any $H$-distinguished representation $\pi$ of $G$ {(i.e.~a representation satisfying $(\pi)^H \neq 0$)} satisfies $\pi \circ \theta \cong \pi^*$.
\item  For any $g\in G$ there are $h_1,h_2\in H$ s.t. $h_1g^{-1} h_2 =\theta(g)$
\end{enumerate}
\end{theorem}
This theorem implies the previous one.

We can generalize this Theorem  further:

\begin{theorem}
\label{thm:GK3}
Let $G$ be a finite group,  $\theta:G\to G$ be an involution and $H \subset G$ be a subgroup.
Then the following are equivalent:
\begin{enumerate}
\item $(G,H)$ is a Gelfand pair and any $H$-distinguished representation possesses a symmetric non-zero bilinear form $B$ satisfying  $B(\pi(g)v,w)=B(v,\pi(\theta(g^{-1}))w)$.
\item  For any $g\in G$ there are $h_1,h_2\in H$ s.t. $h_1g^{-1} \theta(h_2) =\theta(g)$.
\end{enumerate}
\end{theorem}

Theorem \ref{thm:GK3} implies Theorem \ref{thm:GK2} by \cite[Proposition 3]{Vin06}.

The last theorem follows from  Proposition \ref{prop:decomposition} below which is a reinterpretation of the original proof of the Gelfand-Kazhdan criterion.
In order to formulate this proposition, we  need to recall the definition of the twisted Frobenius–Schur indicator.
\begin{notation}
Let $\pi \in \irr (G)$  and let $\theta:G\to G$ be an involution. We denote $$\eps_\theta(\pi)=
\begin{cases}
0,& \pi\not\simeq \pi^* \circ \theta, \\
1,  & \pi  \text{ possesses a non-zero symmetric bilinear form $B$ satisfying}\\& B(\pi(g)v,w)=B(v,\pi(\theta(g^{-1}))w), \\
-1,  & \pi  \text{ possesses a non-zero anti-symmetric bilinear form $B$ satisfying }\\& B(\pi(g)v,w)=B(v,\pi(\theta(g^{-1}))w). \\
\end{cases}$$
\end{notation}
\begin{proposition}
\label{prop:decomposition}
Let $G$ be a finite group,  $\theta:G\to G$ be an involution and $H \subset G$ be a subgroup.
Let $V$ be the space of functions on $G$ which are left invariant w.r.t. $H$, right invariant w.r.t. $\theta(H)$  and anti-invariant w.r.t. $\sigma:=\theta \circ \mathrm{inv}$, where $\mathrm{inv}:G\to G$ is the inversion. Then
$$V \cong \left(
\bigoplus_{\eps_\theta(\pi)=0} \pi^H \otimes  (\pi^*)^{\theta(H)}
\right )^{S_2,\mathrm{sign}}
\oplus
\left (
\bigoplus_{\eps_\theta(\pi)=1} \Lambda^2(\pi^H)
\right )
\oplus
\left (
\bigoplus_{\eps_\theta(\pi)=-1}Sym^2(\pi^H)
\right ),
$$
where the action of $S_2$  on $\bigoplus_{\eps_\theta(\pi)=0} \pi^H \otimes  (\pi^*)^{\theta(H)}$ is given by the involution $s(v\otimes w)\mapsto w\otimes v$ where $v\otimes w\in \pi^H\otimes (\pi^*)^{\theta(H)}$  and $w\otimes v\in   (\pi^*)^{\theta(H) } \otimes    \pi^H \cong   (\pi^*\circ \theta)^H  \otimes   (\pi^*\circ \theta)$.
\end{proposition}

\begin{cor}
Using the notations above we have,
$${\#\{O\in H \backslash  G/\theta(H):\sigma(O)\neq O\}}=\sum_{\pi\in \irr(G)} \dim(\pi^H) \left ( \dim(\pi^{\theta(H)})-\eps_\theta(\pi) \right ).$$
\end{cor}

\subsection{The case of $p$-adic and real groups}
Some of the arguments above work also for $l$-groups and even real reductive groups.
First of all, as in the original Gelfand-Kazhdan criterion, the second condition in all three theorems should be replaced by a condition on distributions. Similarly, the space $V$ as above should be replaced by a space of distributions.

The proof of \cite{Vin06} works also for $l$-groups (see \cite[Lemma 3]{Vin06}), and the same argument seems to work for real reductive groups.
Thus, the main difference is in Proposition \ref{prop:decomposition}.   Proposition \ref{prop:decomposition} does not work as is in those cases. However, the construction of the spherical (a.k.a.~relative) character gives an embedding $$ \nu:\left(
\bigoplus_{\eps_\theta(\pi)=0} (\pi^*)^H \otimes  (\tilde \pi^*)^{\theta(H)}
\right )^{S_2,\mathrm{sign}}
\oplus
\left (
\bigoplus_{\eps_\theta(\pi)=1} \Lambda^2((\pi^*)^H)
\right )
\oplus
\left(
\bigoplus_{\eps_\theta(\pi)=-1}Sym^2((\pi^*)^H)
\right)\to V.$$
Using this, the implication $(2) \Rightarrow (1)$
 of Theorems \ref{thm:GK1}, \ref{thm:GK2} and \ref{thm:GK3} follows.
 In fact, this is a reformulation of the classical proof of the Gelfand-Kazhdan criterion, and its extension that was proven in  \cite{JR}. It is reasonable to expect that in many cases $\nu$ has dense image. If this is the case, then Theorems \ref{thm:GK1}, \ref{thm:GK2} and \ref{thm:GK3} hold in the $p$-adic and real settings. Namely,  consider the following:

\begin{defn}
Let $H_1,H_2 \subset G$ be subgroups of an $l$-group or of a real reductive group and let $\Sc^*(G)$ denote the space of Schwartz distributions on $G$. We say that $(G,H_1,H_2)$  satisfies spectral density if the space spanned by spherical characters of irreducible  (admissible) representations of $G$  w.r.t. $H_1,H_2$ is dense in $\Sc^*(G)^{H_1\times H_2}$.
\end{defn}
We prove that a weaker property is satisfied in the $p$-adic case for many cases in \cite[Theorems C and D]{AGS_Z}.

The argument above show the following:
\begin{theorem}
If $G$ is an $l$-group (or real reductive group)  and $H$ is a subgroup s.t.  $(G,H,H)$  satisfies  spectral density, then Theorems \ref{thm:GK1}, \ref{thm:GK2} and \ref{thm:GK3} hold for $G,H$ with the above mentioned changes.
\end{theorem}

\section{Double invariant functions and multiplicities (proof of Proposition \ref{prop:decomposition})}
Let $X=G/H$, and set $H'=\theta(H)$ and $X'=G/H'$.
Let $\sigma$ be the involution of $X\times X'$ given by $([g],[h])\mapsto ([\theta(h)],[\theta(g)])$.
We have $$\C[X]=\bigoplus_{\pi\in \irr (G)} \pi \otimes (\pi^*)^H ~~~~~\text{~ and ~}~~~~~ \C[X']=\bigoplus_{\pi\in \irr (G)} \pi \otimes (\pi^*)^{H'}.$$ Thus

$$W:=\C[G]^{H \times H'}\cong \C[X\times X']^{\Delta G}\cong \bigoplus_{\pi \in \irr (G)}  (\pi^*)^H \otimes (\pi)^{H'},$$
where $\Delta G$ denotes the diagonal embedding of $G$ into $G \times G$.
Its remains to understand the action of $\sigma$ on $W$.
For this let us first analyze the action of $\sigma$ on $\C[X\times X']$.   We have

\begin{flalign*}
\C[X\times X']
&\cong \bigoplus_{\pi,\tau \in \irr (G)} \pi\otimes (\pi^*)^H \otimes \tau \otimes (\tau^*)^{H'}
\cong \bigoplus_{\pi,\tau \in \irr (G)} (\pi \circ \theta) \otimes (\pi^*)^{H'} \otimes \tau \otimes (\tau^*)^{H'} \\
&\cong \bigoplus_{\pi,\tau \in \irr (G)} (\pi \circ \theta) \otimes \tau \otimes (\pi^*)^{H'}  \otimes (\tau^*)^{H'}\\
&\cong\left(\bigoplus_{\pi \in \irr (G)} (\pi \circ \theta) \otimes \pi \otimes (\pi^*)^{H'}  \otimes (\pi^*)^{H'}\right)\oplus  \left(\bigoplus_{\pi\not\simeq \tau \in \irr (G)} (\pi \circ \theta) \otimes \tau \otimes (\pi^*)^{H'}  \otimes (\tau^*)^{H'}\right).
\end{flalign*}

The action of $\sigma$ on $$\bigoplus_{\pi\not\simeq \tau \in \irr (G)} (\pi \circ \theta) \otimes \tau \otimes (\pi^*)^{H'}  \otimes (\tau^*)^{H'}$$ is given by interchanging the summand corresponding to $(\pi,\tau)$ with the summand corresponding to $(\tau,\pi)$.
The action of $\sigma$ on
$$\bigoplus_{\pi\in \irr (G)} (\pi \circ \theta) \otimes \pi \otimes (\pi^*)^{H'}  \otimes (\pi^*)^{H'}$$ is by acting on each  summand separately, and is given by $$v\otimes w\otimes \alpha\otimes \beta\mapsto w\otimes v\otimes  \beta \otimes  \alpha.$$ Now, let us restrict  this action to  $((\pi \circ \theta) \otimes \pi \otimes (\pi^*)^{H'}  \otimes (\pi^*)^{H'})^{\Delta G}\cong((\pi \circ \theta) \otimes \pi)^{\Delta G} \otimes (\pi^*)^{H'}  \otimes (\pi^*)^{H'}$.  We see that the space $(\pi \circ \theta) \otimes \pi)^G$ is either $0$ or $1$-dimensional, and in the latter case, the action of $\sigma$ on it is given by $\eps_\theta(\pi)$.

Finally we use the fact that $V=\{ w\in W:\sigma(w)=-w\}$ and obtain the required identity.



\bibliographystyle{alpha}
\bibliography{Ramibib}
\end{document}